%% file: fourier.tex
\documentclass[11pt]{article}
\input{format} \input{mathdefs} \input{thmdefs}

\title{Amenability and weak amenability of the Fourier algebra}
\author{{\it Brian E.\ Forrest}\thanks{Research supported by NSERC under grant no.\ 90749-00.}
\and
{\it Volker Runde}\thanks{Research supported by NSERC under grant no.\ 227043-00.}}
\date{}
\begin{document}
\maketitle
\begin{abstract}
Let $G$ be a locally compact group. We show that its Fourier algebra $A(G)$ is amenable if and only if $G$ has an
abelian subgroup of finite index, and that its Fourier--Stieltjes algebra $B(G)$ is amenable if and only if $G$ has
a compact, abelian subgroup of finite index. We then show that $A(G)$ is weakly amenable if the component of the identity of $G$ is
abelian, and we prove some partial results towards the converse.
\end{abstract}
\section*{Introduction}
In \cite{Joh1}, B.\ E.\ Johnson initiated the theory of amenable
Banach algebras. He proved that the amenable locally compact groups 
$G$ can be characterized through a cohomological triviality condition for the
group algebra $L^1(G)$. This triviality condition can be formulated
for any Banach algebra and is used to define the class of amenable Banach algebras.
\par
The Fourier algebra $A(G)$ of a locally compact group $G$ was
introduced by P.\ Eymard in \cite{Eym}. 
Since an amenable Banach algebra always has a bounded approximate identity, Leptin's theorem (\cite{Lep})
immediately yields that $A(G)$ can be amenable only for amenable $G$. 
\par
The first to systematically investigate which locally compact groups
have amenable Fourier algebras was Johnson in
\cite{Joh3}. For any locally compact group $G$, let $\hat{G}$ denote
its dual object, i.e.\ the collection of all equivalence classes of
irreducible unitary representations of $G$. For $\pi \in \hat{G}$, let
$d_\pi$ denote its degree, i.e.\ the dimension of the corresponding
Hilbert space. For compact $G$, Johnson showed: If $\sup \{ d_\pi :
\pi \in \hat{G} \} < \infty$, then $A(G)$ is amenable, whereas, for infinite $G$
such that $\{ \pi \in \hat{G} : d_\pi = n \}$ is finite for each $n
\in \posints$, the Fourier algebra cannot be amenable. Hence, for
example, $A(\operatorname{SO}(3))$ is not amenable; in fact, it is not
even weakly amenable (\cite[Corollary 7.3]{Joh3}).
\par
Soon after Johnson, A.\ T.-M.\ Lau, R.\ J.\ Loy, and G.\ A.\ Willis
(\cite{LLW}) --- and, independently, the first-named author in
\cite{For} --- extended one of Johnson's results by showing that any
locally compact group $G$ having an abelian subgroup of finite index has
an amenable Fourier algebra (this condition is equivalent to $\sup \{ d_\pi :
\pi \in \hat{G} \} < \infty$; see \cite{Moo}).
\par
As the predual of the group von Neumann algebra $\VN(G)$, the Fourier
algebra $A(G)$ has a canonical operator space structure with respect
to which multiplication is completely contractive (see
\cite{ER}). Adding operator space overtones to Johnson's definition of
an amenable Banach algebra, Z.-J.\ Ruan (\cite{Rua}) introduced a variant of
amenability --- called operator amenability --- for which an analog of
Johnson's theorem for $A(G)$ is true: A locally compact group $G$ is
amenable if and only if $A(G)$ is operator amenable (\cite[Theorem 3.6]{Rua}).
There are further results that suggest that in order to deal with
amenability and its variants for the Fourier algebra, one has to take
the canonical operator space structure into account. For example, N.\
Spronk proved in \cite{Spr} that $A(G)$ is operator weakly amenable
for every locally compact group $G$.
\par
Nevertheless, the question for which locally compact groups $G$
precisely the Fourier algebra $A(G)$ is amenable (\cite[Problem 14]{LoA}) or weakly amenable remains an intriguing open problem. No
locally compact groups  besides those with an abelian subgroup of
finite index are known to have an amenable Fourier algebra. In view of
\cite{Joh3} and \cite{Los}, it is plausible to conjecture that these
are the only ones. In the present paper, we prove this conjecture. 
\par
As far as the weak amenability of the Fourier algebra is concerned, we prove that, if the component of the identity of a locally compact grup $G$ is abelian, then
$A(G)$ is weakly amenable. We believe, but have been unable to prove, that the converse holds as well. We introduce the formally stronger notion of hereditary weak
amenability, and then show for $[\SIN]$-groups, that the hereditary weak amenability of $A(G)$ forces the component of the identity of $G$ to be abelian.
\par
Part of the material in this paper has been taken from the unpublished
manuscripts \cite{For} and \cite{RunNon} which it supersedes.
\section{Preliminaries}
Let $\A$ be a Banach algebra, and let $E$ be a Banach $\A$-bimodule. A
{\it derivation\/} $D \!: \A \to E$ is a bounded, linear map
satisfying
\[
  D(ab) = a \cdot Db + (Da) \cdot b \qquad (a,b \in \A).
\]
A derivation $D$ is called {\it inner\/} if there is $x \in E$ such
that
\[
  Da = a \cdot x - x \cdot a \qquad (a \in \A).
\]
The dual space $E^\ast$ of $E$ becomes a Banach $\A$-bimodule via
\[
  \langle x, a \cdot \phi \rangle := \langle x \cdot a, \phi \rangle
  \quad\text{and}\quad
  \langle x, \phi \cdot a \rangle := \langle a \cdot x, \phi \rangle
  \qquad (a \in \A, \, \phi \in E^\ast, \, x \in E).
\]
Modules of this kind are called dual Banach $\A$-bimodules.
\par
We recall the definition of an amenable Banach algebra from (\cite{Joh1}):
\begin{definition} \label{amdef}
A Banach algebra $\A$ is said to be {\it amenable\/} if every 
derivation from $\A$ into a dual Banach $\A$-bimodule is inner.
\end{definition}
\par
For the theory of amenable Banach algebra, see \cite{Joh1} and \cite{LoA}.
\par
For many purposes, it is convenient and even necessary not to use
Definition \ref{amdef}, but equivalent, more intrinsic
characterizations. The following is also due to Johnson.
\par
Let $\Tensor$ denote the projective tensor product of Banach
spaces. If $\A$ is a Banach algebra, $\A \Tensor \A$ becomes a Banach $\A$-bimodule via
\[
  a \cdot (x \tensor y) := ax \tensor y \quad\text{and}\quad (x \tensor y) \cdot a := x \tensor ya \qquad (a,x,y \in \A).
\]
Multiplication induces a bounded linear map $\Delta \!: \A \Tensor \A \to \A$. An {\it approximate diagonal\/} for $\A$ is a bounded net $( d_\alpha )_\alpha$ in $\A \Tensor \A$ such that
\[
  a \cdot d_\alpha - d_\alpha \cdot a \to 0 \qquad (a \in \A)
\]
and
\[
  a \Delta d_\alpha \to a \qquad (a \in \A).
\]
The existence of an approximate diagonal characterizes the amenable Banach algebras (\cite{Joh2}):
\begin{theorem} \label{BEJ}
The following are equivalent for a Banach algebra $\A$:
\begin{items}
\item $\A$ is amenable.
\item $\A$ has an approximate diagonal.
\end{items}
\end{theorem}
\par
This characterization allows to introduce a quantitative aspect into
the notion of an amenable Banach algebra: The {\it amenability
  constant\/} of a Banach algebra $\A$ is the smallest $C \geq 1$ such
that $\A$ has an approximate diagonal bounded by $C$.
\par
We also require another characterization of amenable Banach
algebras. For a Banach algebra $\A$, let $\A^\op$ denote the same
algebra with reversed multiplication. Then $\ker \Delta$ becomes a
left ideal in $\A \Tensor \A^\op$.
\par
The following is \cite[Proposition VII.2.15]{Hel}:
\begin{proposition} \label{sasha}
The following are equivalent for a Banach algebra $\A$:
\begin{items}
\item $\A$ is amenable.
\item $\A$ has a bounded approximate identity, and the left ideal
  $\ker \Delta$ of $\A \Tensor \A^\op$ has a bounded right approximate identity.
\end{items}
\end{proposition}
\par
Let $\A$ be a Banach algebra. A Banach $\A$-bimodule $E$ is called {\it symmetric\/} if
\[
  a \cdot x = x \cdot a \qquad (a \in \A, \, x \in E).
\]
\par
The following definition is from \cite{BCD}:
\begin{definition} \label{WAdef}
A commutative Banach algebra $\A$ is said to be {\it weakly amenable\/} if there is no non-zero derivation from $\A$ into a symmetric Banach $\A$-bimodule.
\end{definition}
\par
It is easy to see that weak amenability is indeed weaker than amenability.
\par
Definition \ref{WAdef} is of little use for non-commutative $\A$. Even though it is possible to extend the notion of weak amenability meaningfully to non-commutative Banach algebras,
we content ourselves with this definition: all the algebras we consider are commutative. The hereditary properties of weak amenability for commutative Banach algebras --- which parallel those of
amenability --- are discussed in \cite{Groe}.
\par
The algebras we are concerned with in this paper were introduced in \cite{Eym}. Let $G$ be a locally compact group. Then the {\it Fourier algebra\/} $A(G)$ of $G$ is the collection of all functions
\[
  G \to \comps, \quad x \mapsto \langle \lambda(x) \xi, \eta \rangle,
\]
where $\xi, \eta \in L^2(G)$ and $\lambda$ is the left regular representation of $G$ on $L^2(G)$. The {\it Fourier--Stieltjes algebra\/} $B(G)$ of $G$ consists of all functions
\[
  G \to \comps, \quad x \mapsto \langle \pi(x) \xi, \eta \rangle,
\]
where $\pi$ is some unitary representation of $G$ on a Hilbert space
$\mathfrak H$ --- always presumed to be continuous with respect to the given topology
on $G$ and the weak operator topology on ${\cal B}(\Hilbert)$ --- and $\xi, \eta \in \mathfrak H$. For more information on $A(G)$ and $B(G)$, in particular for proofs that they are indeed Banach
algebras, we refer to \cite{Eym}.
\par
Let $\VN(G)$ denote the {\it group von Neumann algebra\/} of the locally compact group $G$, i.e.\ the von Neumann algebra acting on $L^2(G)$ that is generated by $\lambda(G)$. As the unique predual of $\VN(G)$, the
Fourier algebra $A(G)$ carries a canonical operator space structure (see \cite{ER} for the necessary background from operator space theory; for more on the operator space $A(G)$, see \cite{FW}). Similarly, 
$B(G)$ --- as the dual space of the full group $\cstar$-algebra $\cstar(G)$ --- is an operator space in a canonical manner.
\par
For any function $f$ on a group $G$, we define a function $\check{f}$ by letting $\check{f}(x) := f(x^{-1})$ for $x \in G$. The following is a compilation of known facts which show that the canonical operator space 
structure of the Fourier algebra is considerably finer than the mere Banach space structure:
\begin{proposition} \label{checkprop}
Let $G$ be a locally compact group. Then the algebra homomorphism
\[
  \theta_\ast \!: A(G) \to A(G), \qquad f \mapsto \check{f}
\]
\begin{items}
\item is an isometry,
\item but is completely bounded if and only if $G$ has an abelian subgroup of finite index.
\end{items}
\end{proposition}
\begin{proof}
(i) is well known (see \cite{Eym}), and the ``if'' part of (ii) follows from \cite[Theorem 4.5]{FW}.
\par
For any Banach space $E$ over $\comps$, let $\overline{E}$ denote its complex conjugate space, i.e.\ the complex Banach space for which scalar multiplication is defined via
$\lambda \odot x := \bar{\lambda}x$ for $\lambda \in \comps$ and $x \in E$. It is easy to see that pointwise conjugation is a complete isometry from $A(G)$ to $\overline{A(G)}$. Consequently, 
if $\theta_\ast$ is completely bounded, so is 
\[
  \bar{\theta}_\ast \!: A(G) \to \overline{A(G)}, \quad f \mapsto \bar{\check{f}}.
\]
with $\| \bar{\theta}_\ast \|_{\mathrm{cb}} = \| \theta_\ast \|_{\mathrm{cb}}$.
The adjoint of $\bar{\theta}_\ast$, however, is nothing but forming the Hilbert space adjoint
\[
  \bar{\theta} \!: \VN(G) \to \overline{\VN(G)}, \quad T \mapsto T^\ast,
\]
which then also has to be completely bounded (with the same $\mathrm{cb}$-norm) by \cite[Proposition 3.2.2]{ER}. 
\par
We claim that this forces $\VN(G)$ to be subhomogeneous, i.e.\ the irreducible representations of $\VN(G)$ are finite-dimensional, and their degrees are bounded. Assume otherwise.
The structure theory of von Neumann algebras then entails that $\VN(G)$ contains the full matrix algebra $\mathbb{M}_n$ as a $^\ast$-subalgebra for each 
$n \in \posints$ (see, e.g., \cite[Chapter 6]{LoA} for details). For $n \in \posints$, let $\theta_n \!: {\mathbb M}_n \to {\mathbb M}_n$ stand for taking the transpose of an $n \times n$-matrix. 
Since entrywise conjugation of matrices is a complete isometry, it follows from \cite[Proposition 2.2.7]{ER} that
\[
  n = \| \theta_n \|_{\mathrm{cb}} = \| \bar{\theta} |_{\mathbb{M}_n} \|_{\mathrm{cb}} \leq \| \bar{\theta} \|_{\mathrm{cb}} = \| \bar{\theta}_\ast \|_{\mathrm{cb}} = 
      \| \theta_\ast \|_{\mathrm{cb}} \qquad (n \in \posints),
\]
which is impossible. Hence, $\VN(G)$ is subhomogeneous.
\par
Let $m \in \posints$ be such that every irreducible representation of $\VN(G)$ has degree $m$ or less. Then $\VN(G)$ satisfies the (non-commutative) polynomial identity $S_{2m} = 0$ (see \cite{AD}
for details), as does its subalgebra $L^1(G)$. By \cite[Lemma 3.4(i)]{AD}, the full group $\cstar$-algebra $\cstar(G)$ then also satisfies $S_{2m} = 0$, which, by \cite[Lemma 3.4(ii)]{AD}, means that
the degrees of the irreducible representations of $\cstar(G)$ are bounded by $m$ as well. This, in turn, entails that $\sup \{ d_\pi : \pi \in \hat{G} \} \leq m$, so that $G$ must have an abelian subgroup 
of finite index by \cite[Theorem 1]{Moo}.
\par
This proves the ``only if'' part of (ii).
\end{proof}
\section{(Non-)amenability of $A(G)$ and $B(G)$}
To prove that every locally compact group with amenable Fourier algebra has an abelian subgroup of finite index, we proceed rather indirectly. 
\par
The {\it anti-diagonal\/} of a (discrete) group $G$ is defined as
\[
  \Gamma_G := \{ (x, x^{-1}) : x \in G \}.
\]
It is clear that $\Gamma_G$ is a subgroup of $G \times G$ if and only
if $G$ is abelian. 
\par
For certain, possibly non-abelian $G$, the anti-diagonal is at least
not very far from being a subgroup:
\begin{example}
Let $G$ be a group with an abelian subgroup, say $A$, of finite
index. Let $x_1, \ldots, x_n \in G$ be representatives of the left
cosets of $A$. For $j =1, \ldots, n$, define
\[
  A_j := \left\{ (a, x_j a^{-1} x^{-1}_j) : a \in A \right\},
\]
so that
\[
  \Gamma_G = (x_1,x_1^{-1}) A_1 \cup \cdots \cup (x_n,x_n^{-1}) A_n.
\]
Hence, $\Gamma_G$ is at least a finite union of left cosets in $G\times G$.
\end{example}
\par
The {\it coset ring\/} $\Omega(G)$ of a group $G$ is the ring of
subsets generated by all left cosets of subgroups of $G$. The example
shows that $\Gamma_G \in \Omega(G \times G)$ if $G$ has an abelian
subgroup of finite index. Our first goal in this section, is to prove the converse, namely that $\Gamma_G$ lies in the coset ring of
$G \times G$ {\it only\/} if $G$ has such a subgroup.
\par
The coset ring of a group is crucial in the definition of piecewise affine maps (see, \cite{Hos}, for example): Let $G$ and $H$ be groups, and let $S \in \Omega(H)$. Then a map $\alpha \!: S \to G$ is called 
{\it piecewise affine\/} if:
\begin{alphitems}
\item there are pairwise disjoint $S_1, \ldots, S_n \in \Omega(H)$ with $S = \bigcup_{j=1}^n S_j$;
\item for each $j =1, \ldots, n$, there is a subgroup $H_j$ of $H$ and an element $x_j \in H$ such that $S_j \subset x_j H_j$;
\item for each $j =1, \ldots, n$ and with $x_j$ and $H_j$ as before, there is a group homomorphism $\beta_j \!: H_j \to G$ and an element $y_j \in G$ such that
$\alpha(y) = y_j \beta_j(x^{-1}_j y)$ for $y \in S_j$.
\end{alphitems}
\par
The following lemma is contained in \cite[Proposition 3.1]{IS} (and originally from M.\ Ilie's thesis \cite{Ili}). We quote it here for convenience:
\begin{lemma} \label{monica}
Let $G$ and $H$ be discrete groups, let $S \in \Omega(H)$, let $\alpha \!: S \to G$ be piecewise affine, and let $\theta \!: A(G) \to B(H)$ be given by
\[
  \theta(f)(x) := \left\{ \begin{array}{cl} f(\alpha(x)), & x \in S, \\ 0, & x \notin S, \end{array} \right.
\]
for $f \in A(G)$ and $x \in H$. Then $\theta$ is a completely bounded algebra homomorphism.
\end{lemma}
\par
With the help of this lemma, we can now prove:
\begin{proposition} \label{antid}
The following are equivalent for a (discrete) group $G$:
\begin{items}
\item $\Gamma_G \in \Omega(G \times G)$.
\item $G$ has an abelian subgroup of finite index.
\end{items}
\end{proposition}
\begin{proof}
In view of the example at the beginning of this section, only (i) $\Longrightarrow$ (ii) needs proof.
\par
Suppose that $\Gamma_G \in \Omega(G \times G)$. From \cite[4.3.1, Theorem]{Rud}, it follows that 
\[
  G \to G , \quad x \mapsto x^{-1}
\]
is piecewise affine. (Even though \cite[4.3.1, Theorem]{Rud} is stated and proved for abelian groups only, it is true for arbitrary $G$; see \cite[Lemma 1.2(ii)]{IS}.)
Lemma \ref{monica} then yields that $\theta_\ast \!: A(G) \to A(G)$ in Proposition \ref{checkprop} is completely bounded, so that
$G$ must have an abelian subgroup of finite index by Proposition \ref{checkprop}(ii).
\end{proof}
\begin{remark}
Since the sets in the coset ring of a group allow for a rather explicit description (see \cite{Brian}), a more elementary proof of Proposition \ref{antid} should be
possible. 
\end{remark}
\par
We shall now characterize those locally compact groups with an amenable Fourier algebra:
\begin{theorem} \label{thm}
The following are equivalent for a locally compact group $G$:
\begin{items}
\item $A(G)$ is amenable.
\item $G$ has an abelian subgroup of finite index.
\end{items}
\end{theorem}
\begin{proof}
In view of \cite[Theorem 4.1]{LLW}, only (i) $\Longrightarrow$ (ii) needs
proof.
\par
By Proposition \ref{sasha}, the kernel of $\Delta \!: A(G) \Tensor A(G) \to
A(G)$ has a bounded approximate identity, say $( u_\alpha )_{\alpha \in \mathbb{A}}$. For each $\alpha \in \mathbb A$, there are
sequences $\left( f^{(\alpha)}_k \right)_{k=1}^\infty$ and 
$\left( g^{(\alpha)}_k \right)_{k=1}^\infty$ in $A(G)$ such that
\[
  \sum_{k=1}^\infty \left\|  f^{(\alpha)}_k \right\| \left\|
  g^{(\alpha)}_k \right\| < \infty
  \qquad\text{and}\qquad
  u_\alpha = \sum_{k=1}^\infty f^{(\alpha)}_k \tensor g^{(\alpha)}_k.
\]
Thanks to Proposition \ref{checkprop}(i), we may define a bounded net $( v_\alpha )_{\alpha
  \in \mathbb{A}}$ in $A(G) \Tensor A(G)$ by letting
\[
  v_\alpha = \sum_{k=1}^\infty f^{(\alpha)}_k \tensor
  \check{g}^{(\alpha)}_k \qquad (\alpha \in \mathbb{A}).
\]
It is immediate that $(v_\alpha )_{\alpha \in \mathbb{A}}$ is a bounded approximate identity for the kernel of
\[
  \Gamma \!: A(G) \Tensor A(G) \to A(G), \quad f \tensor g \mapsto f \check{g}.
\]
\par
For $f,g \in A(G)$, the function
\[
  G \times G \to \comps, (x,y) \mapsto f(x) g(y)
\]
lies in $A(G \times G)$: this induces a canonical contraction from $A(G) \Tensor A(G)$ to $A(G \times G)$ (note that, by \cite{Los}, we cannot a priori suppose that
$A(G) \Tensor A(G) \cong A(G \times G)$). Nevertheless, we can use this
canonical inclusion to identify the elements of $A(G) \Tensor A(G)$ with
elements of $A(G \times G)$. Define
\[
  I := \left\{ f \in A(G \times G) : \lim_\alpha f v_\alpha = f \right\}.
\]
Then $I$ is a closed ideal of $A(G \times G)$ which contains (the canonical image of) $\ker \Gamma$ and,
by definition, has a bounded approximate identity. Since $G$ --- and thus
$G \times G$ --- is amenable, \cite[Theorem 2.3]{FKLS} applies. Hence, $I$ must be of the form
\[
  I(F) := \{ f \in A(G \times G) : f |_F \equiv 0 \},
\]
where $F$ is the hull of $I$ in $G \times G$; furthermore, this hull must be contained in $\Omega(G \times G)$. Since $( v_\alpha )_\alpha$ lies in $\ker \Gamma$, 
it is easy to see that $\Gamma_G \subset F$, and since $( v_\alpha )_\alpha$ converges pointwise to the indicator function of
$(G \times G) \setminus \Gamma_G$, the converse inclusion holds as well. Hence, $\Gamma_G = F \in \Omega(G \times G)$ holds. From Proposition \ref{antid}, we conclude that $G$ has an
abelian subgroup of finite index.
\end{proof}
\begin{remarks}
\item The results from \cite{Joh3} indicate that the amenability
  constant of $A(G)$ and the maximum degree of the irreducible
  representations of $G$ are likely to be closely related. It would be interesting
  to further investigate this possible relation.
\item In \cite{RunAp}, the second-named author introduced another
  operator space structure, denoted by $OA(G)$, over the Fourier
  algebra $A(G)$ which, in general, is different from the canonical
  one --- even though both have the same underlying Banach space. The
  operator space $OA(G)$ is a completely contractive Banach algebra in
  the sense of \cite{Rua}. At the end of \cite{RunAp}, the author
  conjectured that $OA(G)$ is operator amenable if and only if $G$ is
  amenable. Thanks to \cite[Proposition 2.4]{RunAp}, however, the proof of
  Theorem \ref{thm} can easily be adapted to yield that $OA(G)$ is
  operator amenable if and only if $G$ has an abelian subgroup of
  finite index.
\end{remarks}
\par
With Theorem \ref{thm} proven, the corresponding assertion for Fourier--Stieltjes algebras is easy to obtain:
\begin{corollary} \label{cor}
The following are equivalent for a locally compact group $G$:
\begin{items}
\item $B(G)$ is amenable.
\item $G$ has a compact, abelian subgroup of finite index.
\end{items}
\end{corollary}
\begin{proof}
(i) $\Longrightarrow$ (ii): Since $A(G)$ is a complemented ideal in $B(G)$, the hereditary properties of amenability
immediately yield the amenability of $A(G)$ (\cite[Theorem 2.3.7]{LoA}), so that $G$ has an abelian subgroup $A$ of finite index by Theorem \ref{thm}. We may replace $A$ by its closure and
thus suppose that $A$ is closed. Since $A$ is then automatically open, the restriction map from $B(G)$ to $B(A)$ is surjective. Hence, $B(A)$ is also amenable by \cite[Proposition 2.3.1]{LoA}.
Let $\hat{A}$ denote the dual group of $A$. The Fourier--Stieltjes
transform yields an isometric algebra isomorphism of $B(A)$ and the measure
algebra $M(\hat{A})$, so that $M(\hat{A})$ is amenable. Amenability of $M(\hat{A})$, however, forces $\hat{A}$ to be
discrete (\cite{BM}; see also \cite{DGH} for a more general result). Hence, $A$ must be compact.
\par
(ii) $\Longrightarrow$ (i): If $G$ has a compact, abelian subgroup of finite index, then $G$ itself is compact, so that $B(G) = A(G)$ is amenable by Theorem \ref{thm}.
\end{proof}
\begin{remark}
Corollary \ref{cor} characterizes the locally compact groups $G$ for which $B(G)$ is amenable. The corresponding --- and in a certain sense more natural --- characterization for operator amenability seems to
by still open. A partial result can be found in \cite{RS}.
\end{remark}
\par
Another consequence of Theorem \ref{thm} --- slightly extending it --- is:
\begin{corollary}
The following are equivalent for a locally compact group $G$ with at least two elements:
\begin{items}
\item $G$ has an abelian subgroup of finite index.
\item $A(G)$ is amenable.
\item $I(H)$ is amenable for each closed proper subgroup $H$ of $G$.
\item $I( \{ e \})$ is amenable.
\item $I(H)$ is amenable for some closed proper subgroup $H$ of $G$.
\end{items}
\end{corollary}
\begin{proof}
(i) $\Longleftrightarrow$ (ii) is Theorem \ref{thm}.
\par
(ii) $\Longrightarrow$ (iii): By \cite[Theorem 1.5]{FKLS}, $I(H)$ has a bounded approximate identity and thus is amenable by \cite[Theorem 2.3.7]{LoA}.
\par
(iii) $\Longrightarrow$ (iv) $\Longrightarrow$ (v) are trivial.
\par
(iv) $\Longrightarrow$ (ii): Let $x \in G \setminus H$, so that $xH \cap H = \void$. Then $I(xH) \cong I(H)$ is amenable. Restriction to $H$, maps $I(xH)$ onto a dense subalgebra of $A(H)$, so that $A(H)$
is amenable. Since $I(H)$ is just the kernel of the restriction map, the amenability of $A(G)$ follows from \cite[Theorem 2.3.10]{LoA}. 
\end{proof}
\section{Weak amenability of $A(G)$}
In \cite{Joh3}, Johnson showed that $A(G)$ is weakly amenable whenever $G$ is a totally disconnected, compact group, and soon thereafter, the first-named author extended this result to arbitrary totally disconnected, locally compact
groups (\cite[Theorem 5.3]{cosets}): this shows already that weak amenability of the Fourier algebra imposes considerably weaker constraints on the underlying group than amenability. In this section, we shall present a sufficient condition
--- which we believe to be necessary as well --- which forces a locally compact group to have a weakly amenable Fourier algebra and which subsumes \cite[Theorem 5.3]{cosets}.
\par
We begin with an elementary hereditary lemma:
\begin{lemma} \label{herlem}
Let $G$ be a locally compact group, and let $H$ be a closed subgroup of $G$. Then:
\begin{items}
\item If $A(G)$ is weakly amenable, then so is $A(H)$.
\item If $H$ is open, and if $A(H)$ is weakly amenable, then so is $A(G)$.
\end{items}
\end{lemma}
\begin{proof}
Since $A(H)$ is a quotient of $A(G)$ by \cite[Lemma 3.8]{Brian0}, (i) follows immediately from the basic hereditary properties of weak amenability for commutative Banach algebras (\cite{Groe}).
\par
For (ii), let $E$ be a symmetric Banach $A(G)$-bimodule, and let $D \!: A(G) \to E$ be a derivation. Let $f \in A(G)$ have compact support, and let $x_1, \ldots, x_n \in G$ be such that $\supp(f) \subset
x_1 H \cup \cdots \cup x_n H$. For $j =1, \ldots, n$, let $\chi_j \in B(G)$ be the indicator function of $x_j H$. Then $\A_j := \chi_j A(G)$ is a commutative subalgebra of $A(G)$ isometrically isomorphic to $A(H)$.
Since $A(H)$ is weakly amenable, $D |_{\A_j} \equiv 0$ for $j=1, \ldots, n$ follows. In particular, $D(\chi_j f) = 0$ holds for $j=1, \ldots, n$; since $f = \chi_1 f + \cdots + \chi_n f$, it follows that $Df = 0$.
Since the elements with compact support are dense in $A(G)$, we obtain that $D$ is zero on all of $A(G)$.
\end{proof}
\par
For every locally compact group $G$, we shall denote the component of $G$ containing the identity element by $G_e$. It is well known that $G_e$ is a closed, normal subgroup of $G$. 
Our main result in this section (Theorem \ref{WAthm} below) is that that $A(G)$ is weakly amenable whenever $G_e$ is abelian. 
\par
We first prove a lemma:
\begin{lemma} \label{WAlem}
Let $G$ be a locally compact group such that $G_e$ is abelian, and let $K$ be a compact, normal subgroup of $G$ such that $G/K$ is a Lie group. Then $A(G/K)$ is weakly amenable.
\end{lemma}
\begin{proof}
The component of the identity of $G/K$ is $\cl{\pi(G_e)}$, where $\pi \!: G \to G/K$ is the quotient map (\cite[(7.12) Theorem]{HR}). Thus, if $G_e$ 
abelian, so is $(G/K)_e$; in particular, $A((G/K)_e)$ is weakly amenable. Since $G/K$ is a Lie group, $(G/K)_e$ is open. It follows from Lemma \ref{herlem}(ii) that $A(G/K)$ is weakly amenable.
\end{proof}
\begin{theorem} \label{WAthm}
Let $G$ be a locally compact group such that $G_e$ is abelian. Then $A(G)$ is weakly amenable.
\end{theorem}
\begin{proof}
We shall first treat the case where $G$ is almost connected, i.e.\ where $G / G_e$ is compact.
\par
Let $E$ be a symmetric Banach $A(G)$-bimodule, let $D \!: A(G) \to E$ be a derivation, let $f \in A(G)$ be arbitrary, and let $\epsilon > 0$. We claim that $\| D f \| \leq \epsilon$: since $\epsilon > 0$ is arbitrary, this
is enough to conclude that $Df = 0$.
\par
For $x \in G$, define the right translate $R_x f$ of $f$ by $x$ as $(R_xf)(y) := f(yx)$ for $y \in G$. Since $G \ni x \mapsto R_xf \in A(G)$ is continuous with respect to the norm topology on $A(G)$ (\cite{Eym}),
there is a neighborhood $U$ of the identity in $G$ such that $\| f - R_x f \| < \frac{\epsilon}{1+ \| D \|}$ for all $x \in U$. By \cite[12.2.15 Theorem]{Pal}, there is a compact, normal subgroup $K \subset U$ of $G$ 
such that $G/K$ is a Lie group. Define $P_K \!: A(G) \to A(G)$ by letting
\[
  P_K f := \int_K R_x f \, dx \qquad (f \in A(G)),
\]
where the integral is a Bochner integral with respect to normalized Haar measure on $K$. It follows that
\[
  \| f  - P_K f \| \leq \int_K \| f - R_x f \| \, dx \leq \frac{\epsilon}{1+ \|D \|} \qquad (f \in A(G)).
\]
It is easy to see that $P_K$ is a projection onto $A(G:K)$, the subalgebra of $A(G)$ consisting of those functions that are constant on cosets of $K$; this algebra is isometrically isomorphic to $A(G/K)$ (\cite{Eym}).
Since $A(G:K) \cong A(G/K)$ is weakly amenable by Lemma \ref{WAlem}, we have $D |_{A(G:K)} \equiv 0$ and, consequently, 
\[
  \| Df\| = \|Df - D(P_Kf) \| \leq \| D \| \|f - P_K f \| \leq \epsilon.
\]
This proves our claim.
\par
For arbitrary locally compact $G$, note that $G$ has an open, almost connected subgroup $H$ (this follows, for instance, from \cite[(7.7) Theorem]{HR}). Since $A(G)$ is weakly amenable if and only if this is true for $A(H)$ by 
Lemma \ref{herlem}(ii), the claim for general $G$ follows from the almost connected case. 
\end{proof}
\par
As we already stated at the beginning of this section, we believe that the sufficient condition of Theorem \ref{WAthm} for the weak amenability of the Fourier algebra is also necessary, but we have not been able to prove it. 
If this conjecture is correct, it implies that no non-abelian, connected group can have a weakly amenable Fourier algebra. 
\par
The following proposition is a small step into this converse direction:
\begin{proposition} \label{lieprop}
Let $G$ be a Lie group such that $A(G)$ is weakly amenable. Then every compact
subgroup of $G$ has an abelian subgroup of finite index.
\end{proposition}
\begin{proof}
Let $K$ be a compact subgroup of $G$. If $A(G)$ is weakly amenable, then Lemma \ref{herlem}(i) implies that $A(K)$ is also
weakly amenable. For the same reason, $A(K_e)$ is weakly amenable. But $K_e$ is a compact
connected Lie group. It follows from the remarks following 
\cite[Corollary 7.3]{Joh3} that $K_e$ is abelian. Finally, $K_e$ is open
in the compact Lie group $K$ and is therefore of finite index.
\end{proof}
\par
The equivalence of (i) and (iii) in the next corollary can be already be found in 
\cite[Proposition 4.4]{LLW}:
\begin{corollary}
Let $G$ be an almost connected, semisimple Lie group. Then the following are
equivalent: 
\begin{items}
\item $A(G)$ is amenable. 
\item $G$ is amenable, and $A(G)$ is weakly amenable.
\item $G$ is finite.
\end{items}
\end{corollary}
\begin{proof}
(i) $\Longrightarrow$ (ii) is obvious.
\par
For (ii) $\Longrightarrow$ (iii), simply note that, if $G$ is amenable, then it is already compact (\cite[Theorem 3.8]{Pat}). It follows from Proposition \ref{lieprop}
$G_e$ is both semisimple and abelian.
This means that $G_e$ is trivial and, in turn, that $G$ is finite.
\par
(iii) $\Longrightarrow$ (i) is trivial.
\end{proof}
\par
Let $K$ be a compact normal subgroup of the locally compact group $G$. In the proof of Theorem \ref{WAthm}, we made crucial use of the fact that $A(G)$ contains a closed subalgebra $A(G:K)$ which is
isometrically isomorphic to $A(G/K)$. There ought to be a strong connection between the weak amenability of $A(G)$ and that
of $A(G/K)$: otherwise, it would be extremely unlikely that any reasonable structural characterization exists of those $G$ for which $A(G)$ is weakly amenable. However,
since subalgebras generally do not inherit weak amenability, we cannot conclude (yet) that, if $A(G)$ is weakly
amenable, then so is $A(G/K)$. 
\par
This prompts us to formulate the following definition:
\begin{definition} \label{hWA}
Let $G$ be a locally compact group. We say that $A(G)$ is {\it hereditarily weakly amenable\/} if $A(G/K)$ is weakly amenable for
every compact normal subgroup $K$ of $G$.
\end{definition}
\begin{remark}
If $G_e$ is abelian, then so is $(G/K)_e$ for each compact normal subgroup $K$ of $G$. Hence, the only examples of locally compact groups $G$ for which we positively know that $A(G)$
is weakly amenable, also satisfy the formally stronger Definition \ref{hWA}.
\end{remark}
\par
Recall that a locally compact group $G$ is called a $[\SIN]$-group if it has a basis of neighborhoods of the identity consisting of sets that are invariant under conjugation.
\par
With Definition \ref{hWA} in place of mere weak amenability, we can prove a converse of Theorem \ref{WAthm} for $[\SIN]$-groups:
\begin{theorem} \label{SINthm}
Let $G$ be a $[\SIN]$-group. Then $A(G)$ is hereditarily weakly amenable if
and only if $G_e$ is abelian.
\end{theorem}
\begin{proof}
In view of the remark following Definition \ref{hWA}, it is clear that $A(G)$ is hereditarily weakly amenable whenever $G_e$ is abelian (this does not require $G$ to be a $[\SIN]$-group).
\par
Conversely, suppose that $A(G)$ is hereditarily weakly amenable. Since $G$ is a $[\SIN]$-group, there are a vector group $V$ and a compact group $K$ such that
\[
  G_e \cong V \times K
\]
(see \cite[12.4.48 Theorem]{Pal}). If $K$ is not
abelian, we can find $x,y\in K$ and a neighborhood $U$ in $G$ of the identity such
that the commutator $[x,y]$ does not lie in $U$. Since $G$ is a $[\SIN]$-group, $U$ contains
a compact normal subgroup $K_U$ such that $G/K_U$ is a Lie group (see \cite[12.6.12]{Pal}). Let $\pi_U \!: G\to G/K_U$ be the quotient map.
Then $\pi_U(K)$ is a compact connected Lie group. From the choice of $U$, it is clear that $\pi_U(K)$ is not abelian
in $G/K_U$. We know that $A(G/K_U)$ is weakly amenable because $A(G)$ is hereditarily weakly amenable, and Lemma \ref{herlem}(i) implies that $A(\pi_U(K))$ is also
weakly amenable. Since $\pi_U(K)$ is a non-abelian, compact, connected Lie group, this is impossible by Proposition \ref{lieprop}. We conclude that $K$ is abelian, so that $G_e$ is also abelian.
\end{proof}
\par
Recall that a locally compact group $G$ is called an $[\IN]$-group if the identity of $G$ has a compact neighborhood which is invariant under conjugation.
\par
We conclude this paper with a corollary of Theorem \ref{SINthm}:
\begin{corollary}
Let $G$ be a locally compact $[\IN]$-group such that $A(G)$ is hereditarily weakly amenable. Then there exists a compact normal subgroup $K$ of $G_e$ such that $G_e/K$ is abelian.
\end{corollary}

\begin{proof}
Since $G$ is an $[\IN]$-group, there is a compact normal subgroup $N$ of $G$
such that $G/N$ is a $[\SIN]$-group (see \cite[12.1.31 Theorem]{Pal}). It is easy to see that the hereditary weak amenability of $A(G)$ forces $A(G/N)$ to be hereditarily weakly amenable as well. By Theorem \ref{SINthm}, this means that
$(G/N)_e$ is abelian. Letting $K := G_e \cap N$ and observing that $(G/N)_e \cong G_e/K$ yields the result.
\end{proof}
\renewcommand{\baselinestretch}{1.0}
\dated
\vfill
\renewcommand{\baselinestretch}{1.2}
\begin{tabbing}
{\it Second author's address\/}: \= Department of Mathematical and
Statistical Sciences \kill
{\it First author's address\/}: \> Department of Pure Mathematics \\
                                \> University of Waterloo \\
                                \> Waterloo, Ontario \\
                                \> Canada, N2L 3G1 \\[\medskipamount]
{\it E-mail\/}:                 \> {\tt beforres@math.uwaterloo.ca} \\[\bigskipamount]
{\it Second author's address\/}: \= Department of Mathematical and Statistical Sciences \\
                 \> University of Alberta \\
                 \> Edmonton, Alberta \\
                 \> Canada, T6G 2G1 \\[\medskipamount]
{\it E-mail\/}:  \> {\tt vrunde@ualberta.ca} \\[\medskipamount]
{\it URL\/}: \> {\tt http://www.math.ualberta.ca/$^\sim$runde/}
\end{tabbing}                                                                           
\end{document}

%% file: format.tex
\renewcommand{\baselinestretch}{1.2}
\newcommand{\dated}{\mbox{} \hfill {\small [{\tt \today}]}}

%% file: mathdefs.tex
\usepackage{amsmath,amssymb,amscd}
\newcommand{\pf}[1]{\trivlist \item[\hskip\labelsep\it #1\ ]}
\newcommand{\varpf}[1]{\trivlist \item[\hskip\labelsep\sc #1:]}
\newcommand{\qedbox}{$\rlap{$\sqcap$}\sqcup$}
\newcommand{\qed}{\qquad \qedbox \endtrivlist}
\newcommand{\varqed}{\hfill \rule{0.6em}{0.6em} \endtrivlist}
\newenvironment{proof}{\pf{Proof}}{\qed}

\newenvironment{remark}{\pf{Remark}}{\endtrivlist}
\newenvironment{remarks}{\pf{Remarks} 
   \begin{enumerate}}{\end{enumerate} \endtrivlist}
\newenvironment{example}{\pf{Example}}{\endtrivlist}

\newenvironment{items}{
  \begin{enumerate} 
                    
  }{\end{enumerate}}
\newenvironment{alphitems}{
  \begin{enumerate} 
                    
  }{\end{enumerate}}
\newenvironment{keywords}{\noindent\small {\it Keywords\/}:}{\vskip 4pt}
\newenvironment{classification}{\noindent\small 2000 {\it Mathematics Subject
Classification\/}:}{\vskip 12pt}

\newcommand{\comps}{{\mathbb C}}

\newcommand{\posints}{{\mathbb N}}

\newcommand{\void}{\varnothing}
\newcommand{\tensor}{\otimes}

\newcommand{\Tensor}{\hat{\otimes}}

\newcommand{\cstar}{{C^\ast}}

\newcommand{\A}{{\mathfrak A}}

\newcommand{\Hilbert}{{\mathfrak H}}

\newcommand{\op}{{\mathrm{op}}}

\newcommand{\VN}{\operatorname{VN}}

\newcommand{\SIN}{\operatorname{SIN}}
\newcommand{\IN}{\operatorname{IN}}

\newcommand{\supp}{{\operatorname{supp}}}

\newcommand{\cl}[1]{#1^-}

%% file: thmdefs.tex
\newtheorem{theorem}{Theorem}[section]
\newtheorem{lemma}[theorem]{Lemma}
\newtheorem{corollary}[theorem]{Corollary}
\newtheorem{proposition}[theorem]{Proposition}
\newtheorem{df}[theorem]{Definition}
\newenvironment{definition}{\begin{df} \rm}{\end{df}}